\documentclass[12pt]{amsart}
\usepackage{amsfonts}

\usepackage{amsmath}
\usepackage{graphicx}

\newtheorem{theorem}{Theorem}

\newtheorem{corollary}[theorem]{Corollary}

\newtheorem{definition}[theorem]{Definition}
\newtheorem{example}[theorem]{Example}

\newtheorem{lemma}[theorem]{Lemma}

\newtheorem{proposition}[theorem]{Proposition}
\newtheorem{remark}[theorem]{Remark}

\begin{document}
\title[Shape preserving properties]{Shape preserving properties of generalized Bernstein operators on Extended
Chebyshev spaces}
\author[J. M. Aldaz, O. Kounchev and H. Render]{J. M. Aldaz, O. Kounchev and H. Render}
\address{J. M. Aldaz: Departamento de Matem\'aticas, Universidad Aut\'onoma de
Madrid, Cantoblanco 28049, Madrid, Spain.}
\email{jesus.munarriz@uam.es}
\address{O. Kounchev: Institute of Mathematics and Informatics, Bulgarian Academy of
Sciences, 8 Acad. G. Bonchev Str., 1113 Sofia, Bulgaria.}
\email{kounchev@gmx.de}
\address{H. Render: School of Mathematical Sciences, University College Dublin,
Dublin 4, Ireland.}
\email{render@gmx.de}
\thanks{The first and last authors are partially supported by Grant
MTM2006-13000-C03-03 of the D.G.I. of Spain. The second and third authors
acknowledge support within the project ``Institutes Partnership'' with the
Alexander von Humboldt Foundation (Bonn) and within the Project DO-02-275
with the NSF of Bulgarian Ministry of Education and Science. }

\begin{abstract}
We study the existence and shape preserving properties of a generalized
Bernstein operator $B_{n}$ fixing a strictly positive function $f_{0}$, and
a second function $f_{1}$ such that $f_{1}/f_{0}$ is strictly increasing,
within the framework of extended Chebyshev spaces $U_{n}$. The first main
result gives an inductive criterion for existence: suppose there exists a
Bernstein operator $B_{n}:C[a,b]\rightarrow U_{n}$ with strictly increasing
nodes, fixing $f_{0}, f_{1}\in U_{n}$. If $U_{n}\subset U_{n + 1}$ and $U_{n
+ 1}$ has a non-negative Bernstein basis, then there exists a Bernstein
operator $B_{n+1}:C[a,b]\rightarrow U_{n+1}$ with strictly increasing nodes,
fixing $f_{0}$ and $f_{1}.$ In particular, if $f_{0},f_{1},...,f_{n}$ is a
basis of $U_{n}$ such that the linear span of $f_{0},..,f_{k}$ is an
extended Chebyshev space over $\left[ a,b\right] $ for each $k=0,...,n$,
then there exists a Bernstein operator $B_{n}$ with increasing nodes fixing $%
f_{0}$ and $f_{1}.$ The second main result says that under the above
assumptions the following inequalities hold 
\begin{equation*}
B_{n}f\geq B_{n+1}f\geq f
\end{equation*}
for all $(f_{0},f_{1})$-convex functions $f\in C\left[ a,b\right] .$
Furthermore, $B_{n}f$ is $(f_{0},f_{1})$-convex for all $(f_{0},f_{1})$%
-convex functions $f\in C\left[ a,b\right] .$
\end{abstract}

\maketitle

\section{Introduction}

Given $n\in \mathbb{N}$, the space of polynomials generated by $\{1,x,\dots
,x^{n}\}$ on $[a,b]$ is basic in approximation theory and numerical
analysis, so generalizations and modifications abound. However, from a
numerical point of view it is a well known fact that the Bernstein bases
functions $p_{n,k}=x^{k}\left( 1-x\right) ^{n-k}$ behave much better and,
in the sense of \cite{FaGo96},
provide optimal stability. The associated Bernstein
operator $B_{n}:C\left[ 0,1\right] \rightarrow U_{n},$ defined by 
\begin{equation}
B_{n}f\left( x\right) =\sum_{k=0}^{n}f\left( \frac{k}{n}\right) \binom{n}{k}%
x^{k}\left( 1-x\right) ^{n-k}  \label{defBP}
\end{equation}
has been the object of intensive research. As is well known, the polynomials 
$B_{n}f$ converge to $f$ uniformly although the convergence might be very
slow. More important is the fact that the Bernstein operator $B_{n}$ reduces
the variation and preserves the shape of $f$. In particular, if $f$ is
increasing then $B_{n}f$ is increasing, while if $f$ is convex then $B_{n}f$
is convex, see e.g. \cite{Davis}. And the derivative of $B_{n}f$ of a
function of class $C^{1}$ converges uniformly to $f^{\prime }$, cf. \cite
{Lore86}, pg. 25. For this reason Bernstein bases and operators are
fundamental notions.

In Computer Aided Geometric Design (CADG) one is often interested, for
instance, in rendering circumferences and other shapes not given by
polynomial functions. It is thus natural to try to extend the preceding
theory to more general spaces, containing not only $1,x,\dots ,x^{n}$, but
also, say, sine and cosine functions, while keeping as many of the good
properties of Bernstein bases and operators as possible. If one generalizes
the space of polynomials of degree at most $n$ by retaining the bound on the
number of zeros, one is led to the notion of an \emph{extended Chebyshev
space} (or system) $U_{n}$ of dimension $n+1$ over the interval $\left[ a,b%
\right] $: $U_{n}$ is an $n+1$ dimensional subspace of $C^{n}\left( \left[
a,b\right] \right) $ such that each $f\in U_{n}$ has at most $n$ zeros in $%
\left[ a,b\right] $, counting multiplicities, unless $f$ vanishes
identically. Recently, a rich mathematical literature has emerged concerning 
\emph{generalized Bernstein bases} in the framework of extended Chebyshev
spaces, see \cite{CMP04}, \cite{CMP07}, \cite{Cost00}, \cite{CLM}, \cite
{MPS01}, \cite{Mazu99}, \cite{Mazu04}, \cite{Mazu05}, \cite{Mazu05b}, \cite
{MaPo96}, \cite{Pena02}, \cite{Pena05}, \cite{Zhan96}.

It is well-known that extended Chebyshev spaces possess \emph{non-negative
Bernstein bases}, i.e. collections of non-negative functions $%
p_{n,k},k=0,...,n,$ in $U_{n}$, such that each $p_{n,k}$ has a zero of order 
$k$ at $a$ and a zero of order $n-k$ at $b$, for $k=0,...,n$. Assuming that $%
U_{n}$ has a non-negative Bernstein basis $p_{n,k},k=0,...,n$ over the
interval $\left[ a,b\right] $, it is natural to ask whether one may
associate a Bernstein operator $B_{n}:C\left[ a,b\right] \rightarrow U_{n}$
with properties analogous to the classical operator defined in 
(\ref{defBP}). We consider operators $B_{n}$ of the form 
\begin{equation}
B_{n}\left( f\right) =\sum_{k=0}^{n}f\left( t_{n,k}\right) \alpha
_{n,k}p_{n,k}  \label{eqBern}
\end{equation}
where the nodes $t_{n,0},...,t_{n,n}$ belong to the interval 
$\left[ a,b\right]$, and the weights $\alpha_{n,0},...,\alpha_{n,n}$ are positive. But
it is not obvious how the nodes and weights should be defined. Recall that
the classical Bernstein operator reproduces the constant function $1$ and
the identity function $x$. We mimic this feature by requiring that $B_{n}$
fix two functions $f_{0},f_{1}\in U_{n},$ i.e. that 
\begin{equation}
B_{n}\left( f_{0}\right) =f_{0}\text{ and }B_{n}\left( f_{1}\right) =f_{1},
\label{eqBern2}
\end{equation}
where throughout the paper it is assumed that $f_{0} >0 $ and that 
$f_{1}/f_{0}$ is strictly increasing, unless we explicitly state otherwise.
Functions $f_0$ and $f_1$ satisfying the preceding conditions form a
Haar system, in the terminology of \cite[pg. 25]{Karl68}. Following the Editor's suggestion, we shall call $(f_0, f_1)$ a Haar pair.
We show in Section 2 that after choosing $f_{0}$ and $f_{1}$ in $U_{n}$, the
requirements $B_{n}\left( f_{0}\right) =f_{0}$ and $B_{n}\left( f_{1}\right)
=f_{1}$, if they can be satisfied, uniquely determine the location of the
nodes and the values of the coefficients; in other words, there is at most
one Bernstein operator $B_{n}$ of the form (\ref{eqBern}) satisfying (\ref
{eqBern2}) (observe that the only restriction on the nodes is that they belong to $[a,b]$). 

The question of existence of a Bernstein operator in the above sense is
studied in \cite{AKR07} and \cite{AKR08}. Here we present a new, inductive
criterion for the existence of $B_{n}$, making this paper for the most part
self-contained. Let $f_{0},...,f_{n}\in C^{n}\left[ a,b\right] $ and assume
that for each $k=0,...,n$, the linear space $U_{k} := \left\langle
f_{0},...,f_{k}\right\rangle $, generated by $f_{0},...,f_{k},$ is an
extended Chebyshev space of dimension $k+1$. Then, for every $k=1,...,n$,
there exists a Bernstein operator $B_{k}:C\left[ a,b\right] \rightarrow
U_{k} $ fixing $f_{0}$ and $f_{1}$, whose sequence of nodes is strictly
increasing and interlaces with the nodes of $B_{k-1}$, cf. Corollary \ref
{chain}.

Sections 3 and 4 deal with the shape preserving properties of the
generalized Bernstein operator $B_{n}.$ We shall utilize a generalized
notion of convexity, $(f_{0},f_{1})$-convexity, which, according to \cite
{KaSt}, p. 376, is originally due to Hopf, in 1926, and was later
extensively developed by Popoviciu, especially in the context of Chebyshev
spaces. Ordinary convexity corresponds to $(1,x)$-convexity.

Assume there exists a Bernstein operator $B_{n}:C[a,b]\rightarrow U_{n}$
fixing $f_{0}$ and $f_{1}$. We shall show that if $f\in C\left[ a,b\right] $
is $(f_{0},f_{1})$-convex, then 
\begin{equation*}
B_{n}f\geq f,
\end{equation*}
thus generalizing the same inequality for the standard polynomial Bernstein
operator acting on convex functions. Assume next that $B_{n}$ has strictly
increasing nodes, that $U_{n}\subset U_{n+1}$, and that the latter space has
a non-negative Bernstein basis. From the results in Section 2 we know that
there exists a Bernstein operator $B_{n+1}:C[a,b]\rightarrow U_{n+1}$ fixing 
$f_{0}$ and $f_{1}$. In Section 3 we show that 
\begin{equation*}
B_{n}f\geq B_{n+1}f\geq f
\end{equation*}
for all $(f_{0},f_{1})$-convex functions $f\in C\left[ a,b\right] ,$
generalizing once more the corresponding result for the standard polynomial
Bernstein operator. In Section 4 we prove that under the preceding
hypotheses, $B_{n}$ preserves $\left( f_{0},f_{1}\right) $-convexity, i.e., $%
B_{n}f$ is $(f_{0},f_{1})$-convex for all $(f_{0},f_{1})$-convex functions $%
f\in C\left[ a,b\right] $. A similar result is obtained for the so-called $f_{0}$-monotone functions $f.$ These last results follow from the 
general theory of totally positive bases and their shape preserving
properties.

To put in perspective the inductive existence criterion indicated above,
Section 5 (specially, Theorem \ref{ThmExamp5}) clarifies issues regarding the existence of ``good" Bernstein operators, defined using non-decreasing nodes $t_{n, k}\le t_{n, k + 1}$ in
a suitable interval. We focus on 
the linear space $\mathcal{U}_3$ generated by the functions 
$
1,x,\cos x$, and $\sin x
$
on $[0,b]$. It is well known that normalized, totally positive bases
(such as $\left\{\binom{n}{k}
x^{k}\left( 1-x\right) ^{n-k}\right\}_{k=0}^n$ in the polynomial case) posses
optimality properties from the viewpoint of geometric design and shape preservation (cf.  \cite{CaPe93}, \cite{CaPe94}, \cite
{Good96}). 

It might be though that having such good bases would be enough
to define a Bernstein operator fixing 1 and $x$. But this is not
the case. The space $\mathcal{U}_3$ has  a normalized, totally positive basis for every
$b \in (0, 2\pi)$ (cf. \cite{CMP04}, \cite{MPS01}). However, the existence of 
a Bernstein operator fixing 1 and $x$ imposes the stronger condition $b \in (0, \rho _{0}]$,
where $
\rho _{0}\approx 4.4934  
$
is the first positive zero of $b\mapsto
\sin b-b\cos b$: When $\rho_0 < b < 2\pi$ it is not possible to
find nodes in $[0,b]$ so that $1$ and $x$ are fixed by the operator.
When $b \in (\pi, \rho _{0}]$, the nodes do belong to
$[0, b]$, but they fail to be non-decreasing. Thus, a Bernstein operator
can be defined, but it lacks desirable properties; in particular, it does
not preserve convexity. To ensure the existence of a ``good", convexity preserving Bernstein
operator, the stronger condition $b\le \pi$ must be imposed,
and to have a strictly increasing sequence of nodes we need even more: $b < \pi$.  

This paper is essentially self-contained. For simplicity, we consider only
real valued functions when dealing with existence questions. Regarding shape
preserving properties it is of course natural to consider real-valued rather
than complex-valued functions.

We thank the referees and the Editor, both for their very thorough reading of
this paper and for their many suggestions, which lead to a substantial
rewriting of the present article.

\section{Bernstein operators for Extended Chebyshev Spaces.}

We now introduce the concept of a \emph{Bernstein basis} and of a \emph{%
non-negative Bernstein basis} for a linear subspace $U_{n}\subset C^{n}\left[
a,b\right] $ of dimension $n+1.$ In the literature, the expressions
``Bernstein like basis'' or ``B-basis'' are often used instead of
``Bernstein basis''.

\begin{definition}
\textrm{Let }$U_{n}\subset C^{n}\left[ a,b\right] $\textrm{\ be a linear
subspace of dimension }$n+1.$\textrm{\ \ A Bernstein basis (resp.
non-negative Bernstein basis) for $U_{n}$ is a sequence of functions
(resp. non-negative functions) $p_{n,k},k=0,...,n,$ in $U_{n}$, such that
each $p_{n,k}$ has a zero of exact order $k$ at $a$ and a zero of exact
order $n-k$ at $b$, for $k=0,...,n$. }
\end{definition}

By construction,  $\left(p_{n,k}^{(i)} (a)\right)_{i, k = 0,\ldots, n}$ is a triangular matrix with nonzero
diagonal entries. Hence, 
 a Bernstein basis is indeed a \emph{basis} of the
linear space $U_{n}$. Furthermore, the basis functions are unique up to a
non-zero factor, see e.g. Lemma 19 and Proposition 20 in \cite{Veli07}.

As we indicated in the introduction, extended Chebyshev spaces always have
non-negative Bernstein bases. To make this paper as self-contained as
possible, we briefly indicate the reason: Let $\{h_{0},...,h_{n}\}$ be a
basis for $U_{n}$. To obtain a nonzero function $p_{n,k}$ with (at least) $k$
zeros at $a$ and (at least) $n-k$ zeros at $b$, write 
$p_{n,k}:=a_{0}h_{0}+...+a_{n}h_{n}$. We impose the condition of having $k$
zeros at $a$ (which leads to $k$ equations) and $n-k$ zeros at $b$ (which
gives $n-k$ additional equations). Having $n+1$ variables at our disposal,
there is always a non-trivial solution. The assumption that $U_{n}$ is an
extended Chebyshev space guarantees that $p_{n,k}$ has no more than $n$
zeros, so it has exactly $k$ zeros at $a$ and $n-k$ zeros at $b$. In
particular, $p_{n,k}$ is either strictly positive or strictly negative on $%
(a,b)$. Multiplying by $-1$ if needed, we obtain a non-negative $p_{n,k}$.

In Proposition 3.2 in \cite{Mazu05} it is shown that a subspace $U_{n}\subset C^{n}%
\left[ a,b\right] $ possesses a Bernstein basis $p_{n,k},k=0,...,n$ if and
only if every non-zero $f\in U_{n}$ vanishes at most $n$ times on the \emph{%
set} $\left\{ a,b\right\} $ (and not on the interval $\left[ a,b\right] ).$
We mention that the existence of a Bernstein basis in a space $U_{n}\subset C^{n}\left[ a,b%
\right] $ is a rather weak property; e.g. it does not imply the
non-negativity of the basis functions $p_{n,k},k=0,...,n,$ nor the existence
of Bernstein bases on subintervals $\left[ \alpha,\beta\right] $ of $\left[
a,b\right]$.

The next two results are essential tools and standard techniques in CAGD in
the context of degree elevation.

\begin{proposition}
\label{Prop1}Assume that the linear subspaces $U_{n}\subset U_{n+1}\subset
C^{n+1}\left[ a,b\right] $ possess Bernstein bases $p_{n,k},k=0,...,n$, and $%
p_{n+1,k},k=0,...,n+1$. Then 
\begin{equation}
p_{n,k}=\frac{p_{n,k}^{\left( k\right) }\left( a\right) }{p_{n+1,k}^{\left(
k\right) }\left( a\right) }p_{n+1,k}+\frac{p_{n,k}^{\left( n-k\right)
}\left( b\right) }{p_{n+1,k+1}^{\left( n-k\right) }\left( b\right) }%
p_{n+1,k+1}  \label{eqrecc}
\end{equation}
for each $k=0,...,n.$
\end{proposition}

\begin{proof}
Since $p_{n,k}\in U_{n+1},$ the function $p_{n,k}$ is a linear combination of
the basis functions $p_{n+1,k},k=0,...,n+1.$ Using the fact that $p_{n,k}$ has
exactly $k$ zeros at $a$ and $n-k$ zeros at $b$, we see that $p_{n,k}=\alpha
p_{n+1,k}+\beta p_{n+1,k+1}$ for some $\alpha,\beta\in\mathbb{R}.$ Then
$p_{n,k}^{\left(  k\right)  }=\alpha p_{n+1,k}^{\left(  k\right)  }+\beta
p_{n+1,k+1}^{\left(  k\right)  }$ and inserting $x=a$ yields
\[
\alpha=\frac{p_{n,k}^{\left(  k\right)  }\left(  a\right)  }{p_{n+1,k}%
^{\left(  k\right)  }\left(  a\right)  }.
\]
Similarly, $p_{n,k}^{\left(  n-k\right)  }=\alpha p_{n+1,k}^{\left(
n-k\right)  }+\beta p_{n+1,k+1}^{\left(  n-k\right)  }$ and inserting $x=b$
implies that
\[
\beta=\frac{p_{n,k}^{\left(  n-k\right)  }\left(  b\right)  }{p_{n+1,k+1}%
^{\left(  n-k\right)  }\left(  b\right)  }.
\]
\end{proof}

\begin{lemma}
\label{posend} Under the hypotheses of the preceding proposition, assume
additionally that the functions in the Bernstein bases are non-negative.
Then 
\begin{equation}
\frac{p_{n,k}^{\left( k\right) }\left( a\right) }{p_{n+1,k}^{\left( k\right)
}\left( a\right) } > 0 \mbox{ \ \ \ and \ \ \ }\frac {p_{n,k}^{\left(
n-k\right) }\left( b\right) }{p_{n+1,k+1}^{\left( n-k\right) }\left(
b\right) } > 0  \label{poseq}
\end{equation}
for each $k=0,...,n.$
\end{lemma}

\begin{proof}
If $k = 0$ or $k = n$ the assertion is obvious. If $1\le k \le n$, then the
first inequality in (\ref{poseq}) can be obtained from (\ref{eqrecc}): Divide
both sides by $p_{n+1,k} (x)$, and then let $x\downarrow a$. The second
inequality follows in an analogous way. Alternatively, (\ref{poseq}) can be
derived, without using (\ref{eqrecc}), from the well known and elementary fact
that if $f\in C^{\left(  k\right)  }(I)$ has a zero of order $k $ at $c$,
then
\begin{equation}
k!\cdot\lim_{x\rightarrow c}\frac{f\left(  x\right)  }{\left(  x-c\right)
^{k}}=f^{\left(  k\right)  }\left(  c\right)  . \label{eqLim}%
\end{equation}
Of course, the same formula holds for one side limits.
\end{proof}

Let $U_{n}\subset C^{n}\left[ a,b\right] $ be a linear subspace of dimension 
$n+1$ possessing a non-negative Bernstein basis $p_{n,k},k=0,...,n.$ 
 We now introduce the concept of a Bernstein operator fixing a pair of functions:

\begin{definition}\label{Bopdef}
\textrm{We say that a \emph{Bernstein operator }$B_{n}:C\left[ a,b\right]
\rightarrow U_{n}$ fixing the functions $h$ and $g$ exists, 
if there are points $t_{n,0},...,t_{n,n}\in\left[ a,b\right] $ and
coefficients $\alpha_{n,0},...,\alpha_{n,n}>0$ such that the operator $%
B_{n}:C\left[ a,b\right] \rightarrow U_{n}$ defined by 
\begin{equation}
B_{n}f=\sum_{k=0}^{n}f\left( t_{n,k}\right) \alpha_{n,k}p_{n,k}  \label{defB}
\end{equation}
has the property that 
\begin{equation}
B_{n}h =h \text{ and }B_{n}g=g.  \label{eqfix}
\end{equation}
We say that the sequence of nodes $t_{n,0},...,t_{n,n}\in\left[ a,b\right] $ is 
\emph{ strictly increasing }if } 
\begin{equation*}
t_{n,0}<t_{n,1}<...<t_{n,n}.
\end{equation*}
\end{definition}

While  the strict positivity of the
coefficients $\alpha_{n,k}$ for $k=0,...,n$ is included in our definition
of Bernstein operator,  no
restrictions are imposed on the nodes, save that they belong to $[a,b]$. For a natural
example of a Bernstein operator without strictly increasing nodes, see
Proposition \ref{0j} below, or Theorem \ref{ThmExamp5}, Section 5, in the case $b=\pi$. However, if
the nodes fail to be non-decreasing, then the Bernstein operator
may behave in a pathological way, lacking convexity preserving
properties (cf. Theorem \ref{ThmExamp5}).

Two natural questions arise: when is the existence of a Bernstein operator
guaranteed? and, is the Bernstein operator unique? It turns out that
existence depends on \emph{additional} properties of the space $U_{n}$,
while uniqueness is easy to establish. 

In the development below we always assume that the Bernstein operator fixes a Haar pair $(f_0, f_1)$ (i.e., $f_{0}\in U_{n}$ is strictly positive on $%
\left[ a,b\right] $  and $f_{1}\in U_{n}$ is such
that the function $f_{1}/f_{0}$ is strictly increasing on $[a,b]$). The terms
increasing and decreasing are understood in the non-strict sense. For a
constant $c$ positive means $c > 0$, while for a function $f$ it means $f\ge
0$. Of course, once a Bernstein operator fixes a pair $(f_0, f_1)$, 
it fixes every function in its linear span $\langle f_0, f_1\rangle$.

We will consistently use the following notation. Assume that $p_{j,k},$ $%
k=0,...,j$, is a Bernstein basis of the space $U_{j}$. Given $%
f_{0},f_{1}\in U_{j}$, there exist coefficients $\beta_{j,0},...,\beta_{j,j}$
and $\gamma_{j,0},...,\gamma_{j,j}$ such that 
\begin{equation}
f_{0}\left( x\right) =\sum_{k=0}^{j}\beta_{j,k}p_{j,k}\left( x\right) \text{
and }f_{1}\left( x\right) =\sum_{k=0}^{j}\gamma_{j,k}p_{j,k}\left( x\right) .
\label{eqeq}
\end{equation}
The next lemma answers the question of uniqueness positively. 

\begin{lemma}
\label{nodescoeff} Suppose that the linear subspace $U_{n}\subset C^{n}\left[
a,b\right] $, where $n\geq1$, possesses a non-negative Bernstein basis $%
\{p_{n,k}\}_{k=0}^{n}$.
 If there exists a Bernstein operator $B_{n}:C\left[ a,b\right]
\rightarrow U_{n}$ fixing the functions $f_{0},f_{1}\in U_{n}$, then $%
\beta_{n,k} > 0$ for all $k=0,\dots,n$, where
the coefficients
$\beta_{n,k}$ 
are given by (\ref{eqeq}). Moreover, the nodes of $B_{n}$ are
defined, for $k=0$ and $k=n$, by $t_{n,0}=a$ and $t_{n,n}=b$, and in
general, for $k=0, \dots,n$, by 
\begin{equation}
t_{n,k}:=\left( \frac{f_{1}}{f_{0}}\right) ^{-1}\left( \frac{\gamma_{n,k}}{%
\beta_{n,k}}\right),  \label{nodes}
\end{equation}
with the $\gamma_{n,k}$ given by (\ref{eqeq}).
Furthermore, the coefficients of $B_{n}$ are defined, for $k=0,\dots,n$, by 
\begin{equation}
\alpha_{n,k}:=\frac{\beta_{n,k}}{f_{0}(t_{n,k})}.  \label{coeff}
\end{equation}
In particular, 
\begin{equation}
\alpha_{n,0}=\frac{1}{p_{n,0}(a)}\mbox{ \ \ \ \ and \ \ \ \ }\alpha _{n,n}=%
\frac{1}{p_{n,n}(b)}.  \label{coeff0n}
\end{equation}
\end{lemma}

\begin{proof}
Since $B_{n}(f_{0})=f_{0}$ and
\begin{equation}
f_{0}\left(  x\right)  =\sum_{k=0}^{n}\beta_{n,k}p_{n,k}\left(  x\right)  ,
\label{eqeq0}%
\end{equation}
we have
\[
\sum_{k=0}^{n}f_{0}\left(  t_{n,k}\right)  \alpha_{k}p_{n,k}=\sum_{k=0}%
^{n}\beta_{n,k}p_{n,k}.
\]
This entails that $f_{0}\left(  t_{n, k}\right)  \alpha_{n, k}=\beta_{n,k}$,
since $\{p_{n,k}\}_{k=0}^{n}$ is a basis, and now (\ref{coeff}) follows.
Similarly, from $B_{n}f_{1}=f_{1}$ and
\begin{equation}
f_{1}\left(  x\right)  =\sum_{k=0}^{n}\gamma_{k}p_{n,k}\left(  x\right)
\label{eqeq1}%
\end{equation}
we obtain $f_{1}\left(  t_{n,k}\right)  \alpha_{k}=\gamma_{n,k}.$ Using
$f_{0}>0$ and $\alpha_{n,k}>0$ we see that $\beta_{n,k}>0$. Dividing by
$f_{0}\left(  t_{n,k}\right)  \alpha_{n,k}=\beta_{n,k}$, we find that $t_{n,k}
$ satisfies
\begin{equation}
\frac{f_{1}\left(  t_{n,k}\right)  }{f_{0}\left(  t_{n,k}\right)  }%
=\frac{\gamma_{n,k}}{\beta_{n,k}}, \label{eqtk}%
\end{equation}
and now, since $f_{1}/f_{0}$ is injective, its inverse exists and we get
(\ref{nodes}). Next, inserting $x=a$ in (\ref{eqeq0}) and in (\ref{eqeq1}) we
obtain $f_{0}\left(  a\right)  =\beta_{0}p_{n,0}\left(  a\right)  $ and
$f_{1}\left(  a\right)  =\gamma_{0}p_{n,0}\left(  a\right)  $. Thus
\[
\frac{f_{1}\left(  a\right)  }{f_{0}\left(  a\right)  }=\frac{\gamma_{n,0}%
}{\beta_{n,0}},
\]
and it follows by injectivity that $t_{n,0} = a$. An entirely analogous
argument shows that $t_{n,n} = b$. Since $f_{0}(a) = \beta_{n,0} p_{n,0} (a) =
f_{0} (a) \alpha_{n,0} p_{n,0} (a)$ and $f_{0}(b) = \beta_{n,n} p_{n,n} (b) =
f_{0} (b) \alpha_{n,n} p_{n,n} (b)$, (\ref{coeff0n}) follows.
\end{proof}

Lemma \ref{nodescoeff} tells us that to obtain a Bernstein operator $B_{n}$
fixing $f_{0}$ and $f_{1}$, the nodes $t_{n,k}$ must be the ones given by
equation (\ref{nodes}), and the coefficients $\alpha_{n,k}$ by (\ref{coeff}%
). A simple algebraic manipulation then shows that $B_{n}$ does fix $f_{0}$
and $f_{1}$. To construct $B_{n}$, the difficulty  lies in showing that for $%
k=0,...,n$, the numbers 
\begin{equation*}
\frac{\gamma_{n,k}}{\beta_{n,k}}
\end{equation*}
belong to the image of $[a,b]$ under $f_{1}/f_{0}$, so  the  nodes $t_{n,k}$ can be defined. Even if this is the case, it does not
follow in general that the nodes are increasing (cf. Theorem \ref{ThmExamp5}). It does not seem trivial to characterize the spaces 
$U_{n}\subset C^{n}\left[ a,b\right] $ for which there exist a Bernstein
operator fixing a predetermined Haar pair $f_{0},f_{1}\in U_{n},$ cf. \cite{AKR08}.

Here we present a new, inductive criterion: Existence of $B_n$ with strictly increasing
nodes entails existence of $B_{n+1}$  with strictly increasing nodes.
Furthermore,  the nodes at level $n + 1$ interlace  strictly with
 the nodes at level $n$.

\begin{theorem}
\label{ThmExt} Suppose that the linear subspaces $U_{n}\subset
U_{n+1}\subset C^{n+1}\left[ a,b\right] $, where $n\geq1$, possess
non-negative Bernstein bases $p_{n,k},k=0,...,n$, and $p_{n+1,k},k=0,...,n+1$
respectively. If there exists a Bernstein operator $B_{n}:C\left[ a,b\right]
\rightarrow U_{n}$ fixing the functions $f_{0},f_{1}\in U_{n}$, with
strictly increasing nodes $a=t_{n,0}<t_{n,1}<...<t_{n,n}=b,$ then there
exists a Bernstein operator $B_{n+1}:C\left[ a,b\right] \rightarrow U_{n+1}$
fixing $f_{0},f_{1}$, with strictly increasing and strictly interlacing
nodes ${t}_{n+1,0},...,t_{n+1,n+1}$, that is, 
\begin{equation}
a={t}_{n+1,0}=t_{n,0}<{t}_{n+1,1}<t_{n,1}<{t}_{n+1,2}<t_{n,2}< \cdots <{t}%
_{n+1,n}<t_{n,n}={t}_{n+1,n+1}=b.  \label{interlacing}
\end{equation}
\end{theorem}

\begin{proof}
Let us write $f_{0}=\sum_{k=0}^{n+1}{\beta}_{n+1,k}p_{n+1,k}$ and $f_{1}%
=\sum_{k=0}^{n+1}{\gamma}_{n+1,k}p_{n+1,k}.$ By the preceding lemma, if the
Bernstein operator $B_{n+1}:C\left[  a,b\right]  \rightarrow U_{n+1}$ for
$(f_{0}, f_{1})$ exists, then it has the form
\[
B_{n+1}f:=\sum_{k=0}^{n+1}f\left(  {t_{n+1, k}}\right)  {\alpha}_{n+1,
k}p_{n+1,k},
\]
where the positive coefficients $\alpha_{n+1, k}$ are given by (\ref{coeff})
(with $n+1$ replacing $n$), and the increasing nodes $t_{n+1,k}$ are given by
$t_{n+1,0} = a$, by $t_{n+1,n+1} = b$, and in general, by (\ref{nodes}) when
$k=0,\dots,n + 1$. Thus, we need to show, first, that ${\beta}_{n+1,
0},...,{\beta}_{n+1, n+1} > 0$, in order to get the positivity of the
coefficients $\alpha_{n+1, k}$, and second, that
\begin{equation}
\frac{\gamma_{n, k-1}}{\beta_{n, k-1}}< \frac{{\gamma}_{n+1, k}}{{\beta}_{n+1,
k}}< \frac{\gamma_{n,k}}{\beta_{n,k}}\text{ for }k=1,...,n, \label{eqinter}%
\end{equation}
to obtain the (strict) interlacing property of nodes; note that $\gamma
_{n,0}/\beta_{n,0} = \gamma_{n + 1,0}/\beta_{n+ 1,0}$, since both quantities
equal $f_{1}(a) / f_{0} (a)$, and similarly $\gamma_{n,n}/\beta_{n,n} =
\gamma_{n + 1,n+1}/\beta_{n+ 1,n+1}$ (since both quantities equal $f_{1}(b) /
f_{0} (b)$).

At level $n$, by assumption the Bernstein operator is defined via the
coefficients $\alpha_{n,k} > 0$. Now the argument runs as follows: From the
numbers $\alpha_{n,k}$ we obtain the $\beta_{n,k}$, and from these the
$\beta_{n +1,k}$, which in turn give us the $\alpha_{n + 1,k}$.

Since $\beta_{n,k} = f_{0}(t_{n,k})\alpha_{n,k}$, it follows that $\beta_{n,k}
> 0$. From $f_{0}(a) = \beta_{n+1, 0} p_{n+1, 0} (a)$ and $f_{0}(b) =
\beta_{n+1, n+1} p_{n+1, n+1} (b)$ we see that $\beta_{n+1, 0}> 0$ and
$\beta_{n+1, n+1} > 0$. We show next that for $k=1,\dots,n$,
\begin{equation}
{\beta}_{n+1, k}=\beta_{n, k}\frac{p_{n,k}^{\left(  k\right)  }\left(
a\right)  }{p_{n+1,k}^{\left(  k\right)  }\left(  a\right)  }+\beta_{n,
k-1}\frac{p_{n,k-1}^{\left(  n+1-k\right)  }\left(  b\right)  }{p_{n+1,k}%
^{\left(  n+1-k\right)  }\left(  b\right)  },
\label{betapos} 
\end{equation}
from which the positivity of $\beta_{n+1, 1}, \dots, \beta_{n+1, n}$ follows
by Lemma \ref{posend}. Applying the index raising formula given by Proposition
\ref{Prop1} to $f_{0}=\sum_{k=0}^{n}\beta_{n,k}p_{n,k}$, we see that
\[
f_{0}= \beta_{n,0}\frac{p_{n,k}\left(  a\right)  }{p_{n+1,k}\left(  a\right)
} p_{n+1,0}+\sum_{k=1}^{n}\left[  \beta_{n,k}\frac{p_{n,k}^{\left(  k\right)
} \left(  a\right)  }{p_{n+1,k}^{\left(  k\right)  }\left(  a\right)  }%
+\beta_{n,k-1} \frac{p_{n,k-1}^{\left(  n+1-k\right)  }\left(  b\right)
}{p_{n+1,k}^{\left(  n+1-k\right)  }\left(  b\right)  }\right]  p_{n+1,k}%
\]%
\[
+\beta_{n,n} \frac{p_{n,n}\left(  b\right)  }{p_{n+1,n+1}\left(  b\right)
}p_{n+1,n+1},
\]
and we obtain (\ref{betapos}).

Regarding the interlacing property of nodes, another application of the index
raising formula from Proposition \ref{Prop1}, this time to $f_{1}=\sum
_{k=0}^{n}\gamma_{n,k}p_{n,k}$, yields
\begin{equation}
{\gamma}_{n + 1,k}=\gamma_{n,k}\frac{p_{n,k}^{\left(  k\right)  }\left(
a\right)  }{p_{n+1,k}^{\left(  k\right)  }\left(  a\right)  }+\gamma_{n,
k-1}\frac{p_{n,k-1}^{\left(  n+1-k\right)  }\left(  b\right)  }{p_{n+1,k}%
^{\left(  n+1-k\right)  }\left(  b\right)  }
\label{gammaint} 
\end{equation}
for $k=1,...,n.$ To show that $\frac{\gamma_{n, k-1}}{\beta_{n, k-1}}<\frac{{
\gamma}_{n + 1, k}}{{\beta}_{n+1, k}}$, or equivalently, that $\gamma_{n, k-1}
{\beta}_{n + 1, k}<{\gamma}_{n+1, k}\beta_{n, k-1},$ we use formulas
(\ref{betapos}) and (\ref{gammaint}) to rewrite the latter inequality as
\begin{equation}
\gamma_{n,k-1}\left(  \beta_{n,k}\frac{p_{n,k}^{\left(  k\right)  }\left(
a\right)  }{p_{n+1,k}^{\left(  k\right)  }\left(  a\right)  }+\beta_{n,
k-1}\frac{p_{n,k-1}^{\left(  n+1-k\right)  }\left(  b\right)  }{p_{n+1,k}%
^{\left(  n+1-k\right)  }\left(  b\right)  }\right)  <\left(  \gamma
_{n,k}\frac{p_{n,k}^{\left(  k\right)  }\left(  a\right)  }{p_{n+1,k}^{\left(
k\right)  }\left(  a\right)  }+\gamma_{n,k-1}\frac{p_{n,k-1}^{\left(
n+1-k\right)  }\left(  b\right)  }{p_{n+1,k}^{\left(  n+1-k\right)  }\left(
b\right)  }\right)  \beta_{n,k-1}.
\label{step}
\end{equation}
Simplifying and using $\frac{p_{n,k}^{\left(  k\right)  }\left(  a\right)
}{p_{n+1,k}^{\left(  k\right)  } \left(  a\right)  }>0 $ (by Lemma
\ref{posend}), inequality (\ref{step}) is easily seen to be equivalent to
$\frac{\gamma_{n,k-1}}{\beta_{n,k-1}}< \frac{\beta_{n,k}}{\gamma_{n,k}}$,
which is true by (\ref{nodes}) together with the assumptions that $f_{1}%
/f_{0}$ is increasing and that $t_{n,k-1} < t_{n,k}$.

Inequality $\frac{{\gamma}_{n+1, k}}{{\beta}_{n+1, k}}<\frac{\gamma_{n, k}%
}{\beta_{n,k}}$ is proven in the same way.
\end{proof}

For the next corollary we \emph{do not} a priori assume that $f_{0} > 0$ and 
$f_{1}/f_{0}$ is strictly increasing, since multiplying by $-1$ if needed,
these properties can be obtained from the other assumptions.

\begin{corollary}
\label{chain} Let $f_{0},...,f_{n}\in C^{n}\left[ a,b\right] $, and assume
that the linear spaces $U_{k}$ generated by $f_{0},...,f_{k}$ are extended
Chebyshev spaces of dimension $k+1$ for $k=0,...,n.$ Then for every $%
k=1,...,n$, there exists a Bernstein operator $B_{k}:C\left[ a,b\right]
\rightarrow U_{k}$ fixing $f_{0}$ and $f_{1}$, with strictly increasing
nodes and strictly interlacing with those of $B_{k-1}$.
\end{corollary}

\begin{proof}
Since $U_{0}$ is an extended Chebyshev space over $\left[  a,b\right]  $, the
function $f_{0}$ has no zeros. Multiplying by $-1$ if needed, we may assume
that $f_{0}>0$. Since $U_{1}=\langle f_{0},f_{1}\rangle$ is an extended
Chebyshev space over $\left[  a,b\right]  $ it is easy to see that that
$f_{1}/f_{0}$ is either strictly increasing or strictly decreasing. By
multiplying $f_{1}$ by $-1$ if needed one may assume that $f_{1}/f_{0}$ is
strictly increasing. Let $\{p_{1,0},p_{1,1}\}$ be a non-negative Bernstein
basis for $U_{1}$. We define
\begin{equation}
B_{1}f:=\alpha_{1,0}f\left(  a\right)  p_{1,0}+\alpha_{1,1}f\left(  b\right)
p_{1,1}, \label{dim1}%
\end{equation}
where $\alpha_{1,0}=1/p_{1,0}(a)$ and $\alpha_{1,1}=1/p_{1,1}(b)$. Since both
functions $(B_{1}f_{0}-f_{0})\in U_{1}$ and $(B_{1}f_{1}-f_{1})\in U_{1} $
have a zero at $a$ and another zero at $b$, and $U_{1}$ is an extended
Chebyshev space we see that these functions are zero, so $B_{1}$ fixes $f_{0}$
and $f_{1}$. And now the result follows by inductively applying Theorem
\ref{ThmExt} to each $U_{k+1}$ in the chain $U_{1}\subset U_{2}\subset
...\subset U_{n}.$
\end{proof}

It is well known that given an extended Chebyshev space $U_{n},$ one can
find functions $f_{0},...,f_{n}\in C^{n}\left[ a,b\right] $ such that the
linear spaces $U_{k}$ generated by $f_{0},...,f_{k}$ are extended Chebyshev
spaces of dimension $k+1$ for $k=0,...,n,$ see e.g. Proposition 2.8 of \cite
{Mazu05} (cf. also Definition 2.4 in \cite{Mazu05}). The functions $%
f_{0},...,f_{n}$ can be constructed in the following way:  first one shows
that there exists a strictly larger interval $\left[ a,\beta \right] \supset %
\left[ a,b\right] $ such that $U_{n}$ is an extended Chebyshev system over $%
\left[ a,\beta \right] $ (cf. \cite[p. 351]{Mazu05}). Take now $n$ different
points $\xi _{1},...,\xi _{n}$ in the open interval $\left( b,\beta \right) $%
. For each $k=0,...,n$ define a non-zero function $f_{k}\in U_{n}$ vanishing
 on $\xi _{1},...,\xi _{n-k}$. Then the linear spaces $U_{k}$
generated by $f_{0},...,f_{k}$ are extended Chebyshev spaces over $\left[ a,b%
\right] .$ The disadvantage of this procedure is that the choice of the
functions $f_{0}$ and $f_{1}$ cannot be specified in advance, 
but does depend on the space $U_{n}$.  

Thus, Corollary \ref{chain} implies the next result:

\begin{corollary}
\label{space} Let $n\geq1$ and let $U_{n}$ be an extended Chebyshev space
over $\left[ a,b\right] $. Then it is possible to find a Haar pair $f_{0},
f_{1}\in U_{n}$  and a Bernstein operator $B_{n}:C\left[ a,b\right] \rightarrow
U_{n}$ with strictly increasing nodes, such that $B_{n}$ fixes $f_{0}$ and $%
f_{1}$.
\end{corollary}

\begin{remark}
{\rm Observe that the hypothesis of Corollary \ref{space} is weaker than
that of Corollary \ref{chain}, and so is the conclusion, since $f_{0}$ and $%
f_{1}$ are chosen a posteriori, cf. also the discussion at the beginning of
Section 5. }
\end{remark}

Specializing the preceding results to the case of \emph{%
exponential polynomials}, the conclusions we obtain in the real
case are stronger than those from \cite{AKR07} (however, \cite{AKR07}
deals with the more general complex case).
The space $E_{\left( \lambda _{0},...,\lambda _{n}\right) }$ of \emph{%
exponential polynomials} with \emph{exponents} $\lambda _{0},\dots ,\lambda
_{n}\in \mathbb{C}$ is defined by  
\begin{equation}
E_{\left( \lambda _{0},...,\lambda _{n}\right) }:=\left\{ f\in C^{\infty
}\left( \mathbb{R}\right) :\left( \frac{d}{dx}-\lambda _{0}\right)
\cdots\left( \frac{d}{dx}-\lambda _{n}\right) f=0\right\} .  \label{defexp}
\end{equation}
Exponential polynomials provide natural generalizations of the classical,
trigonometric, and hyperbolic polynomials (see \cite{Schu83}), and the $%
\mathcal{D}$-polynomials considered in \cite{GoNe} and \cite{MoNe00}. They
also furnish (under additional assumptions) examples of extended Chebyshev
spaces recently studied in CAGD for the purpose of representing parametric
curves, see \cite{CMP07}, \cite{CLM}, \cite{Mazu99}. Note in particular that
when $\lambda _{k}=0$ for $0\leq k\leq n$, $\lambda _{n+1}= i$ and $\lambda
_{n+1}=-i$, we have $E_{\left( \lambda _{0},...,\lambda _{n+2}\right)
}=\langle 1,x,\dots ,x^{n},\sin x,\cos x\rangle $, where
$\langle g_1,\dots, g_k\rangle$ denotes the vector space spanned by
$g_1,\dots, g_k$. 

It is well known that $E_{\left( \lambda _{0},...,\lambda _{n}\right) }$ is
an extended Chebyshev space over any compact interval $[a,b]$, for any
choice of $\lambda _{0},...,\lambda _{n}\in \mathbb{R}$ (this is
not true if some of the exponents are complex, but here, we
only consider 
real eigenvalues $\lambda_k$; the reader interested in
complex exponents may want to consult \cite{AKR07}, \cite{AKR08} and \cite
{Veli07}).

The results of the present paper give a simple proof of the existence a Bernstein
operator with strictly  increasing and strictly interlacing nodes, fixing a suitable Haar pair.
 Let us
emphasize that the  interlacing property of nodes does not follow from
\cite{AKR07}.

\begin{theorem}
Let $\lambda _{0},...,\lambda _{n}$ be real numbers, let $a < b$, and for $k=1,\dots,n$, let
$U_k:= E_{\left( \lambda _{0},...,\lambda
_{k}\right)}$ over $[a,b]$. If $\lambda_0 = \lambda_1$, we set 
$f_0(x) := e^{\lambda _{0}x}$ and $f_1 (x) := x e^{\lambda_{0}x}$, 
while if $\lambda_0 \neq \lambda_1$, we set 
$f_0(x) := e^{x \min\{\lambda _{0}, \lambda_1\}}$ and $f_1 (x) := e^{x \max\{\lambda _{0}, \lambda_1\}}$. Then
for each $k=1,...,n$, there is a Bernstein operator $B_{k}:C\left[ a,b\right]
\rightarrow E_{\left( \lambda _{0},...,\lambda _{k}\right) }$ fixing $f_{0}$
and $f_{1}$, such that its sequence of nodes is strictly increasing, and
those at level $k>1$ strictly interlace with the nodes of $B_{k-1}$.
\end{theorem}

\begin{proof} It immediately follows from (\ref{defexp}) that $E_{\left( \lambda _{0},...,\lambda _{n}\right) }$ is invariant under permutations of the eigenvalues
$\lambda_i$, so we can assume, with no loss of 
generality, that
$\lambda_0\le \lambda_1$.  Note that if  $\lambda_0 = \lambda_1$,
then   $U_1 = E_{\left( \lambda _{0}, \lambda
_{0}\right)} = \langle f_0, f_1\rangle$. 
Whether we have $\lambda_0 = \lambda_1$ or $\lambda_0 < \lambda_1$,  $f_0 > 0$ and
$f_1/f_0$ is strictly increasing.
Since the spaces $U_{k}=E_{\left( \lambda _{0},...,\lambda _{k}\right) }$
over $\left[ a,b\right] $  are extended
Chebyshev spaces for each $k=0,...,n,$ the result
immediately follows from Corollary \ref{chain}.
\end{proof}

We finish this section with a proposition illustrating our methods
in the {\em classical} polynomial case. The article 
\cite{Ki} exhibits a sequence of positive linear operators converging to the
identity on $C[0,1]$ and fixing $1$ and $x^{2}$. This sequence is obtained
by replacing $x$ in (\ref{defBP}) with a suitably chosen function $r_{n}(x)$
such that $\lim_{n}r_{n}(x)=x$\textrm{. } It is also possible to fix $1$ and 
$x^{2}$ by using the generalized Bernstein operators considered here. As a
matter of fact, it is possible to fix $f_{0}(x) = 1$ and $f_{1}(x) = x^{j}$
for any $j\geq1$ we wish (of course, if $j=1$ we have the standard case). From
Lemma \ref{nodescoeff} we know how to determine the nodes and the
coefficients, i.e., how $B_{n}$ must be constructed.

On the other hand, we \emph{cannot} use Corollary \ref{chain} to conclude
that such a Bernstein operator $B_{n, 0, j}$ exists (the subscripts 0 and $j$
refer to the exponents of the functions being fixed) since whenever $j > 1$,
the space $U_{1}=\left\langle 1,x^{j}\right\rangle $ is not an extended
Chebyshev space over the \emph{closed} interval $\left[ 0,1\right] $: $x^{j}$
has a zero of order $j$. And unlike the situation considered in Theorem \ref
{ThmExt}, the sequence of nodes we obtain is not strictly increasing: 
given $1<j\leq n$, it is easy to see that $t_{n,0}=\cdots=t_{n,j-1}=0$,
simply by
counting zeros at $a=0$, or by the argument given below.

\begin{proposition}
\label{0j} Fix $j>1$, and let $U_{n}$ be the space of polynomials over $%
[0,1] $ of degree at most $n$. For every $n\geq j$, there exists a Bernstein
operator $B_{n,0,j}:C[0,1]\rightarrow U_{n}$ that fixes $1$ and $x^{j}$, and
converges in the strong operator topology to the identity, as $n\rightarrow
\infty$. The operator $B_{n,0,j}$ is explicitly given by 
\begin{equation*}
B_{n,0,j}f(x)=\sum_{k=0}^{n}f\left( \left(\frac{k(k-1)\cdots(k-j+1)}{%
n(n-1)\cdots(n-j+1)}\right) ^{1/j}\right) \binom{n}{k}x^{k}(1-x)^{n-k}.
\end{equation*}
\end{proposition}

\begin{proof}
For the purposes of this argument we set $p_{n,k}(x):=\binom{n}{k}%
x^{k}(1-x)^{n-k}$ (this differs from the notation used in the introduction for
the classical Bernstein polynomials, but it is more convenient here). The
condition $1=B_{n,0,j}1(x)=\sum_{k=0}^{n}\alpha_{n,k}p_{n,k}$ entails that
$\alpha_{n,k}=1$ for all $n,k$. We use the equality $x^{j}=B_{n,0,j}x^{j}$ to
determine the nodes $t_{n,k}$. Writing 
$$
x^{j}=\sum_{k=0}^{n}\gamma
_{n,k}p_{n,k},
$$
 by (\ref{nodes}) we have $t_{n,k}=\gamma_{n,k}^{1/j}$. 
The coefficients $\gamma_{n,k}$ can be obtained in several ways. A 
rather direct one follows next:
$$
x^{j}
=
x^{j}(x+1-x)^{n-j}
=
x^{j}\sum_{k=0}^{n-j}\binom{n-j}{k}x^{k}(1-x)^{n-j-k}
$$
$$ 
=
\sum_{k=j}^{n}\binom{n-j}{k-j}x^{k}(1-x)^{n-k}
=
\sum_{k=j}^{n}\frac{k(k-1)\cdots
(k-j+1)}{n(n-1)\cdots(n-j+1)}\binom{n}{k}x^{k}(1-x)^{n-k}.
 $$
Thus, 
$\gamma_{n,k}=0$ if 
$0\leq k<j$ and 
\begin{equation}
\gamma_{n,k}=\frac{k(k-1)\cdots(k-j+1)}{n(n-1)\cdots(n-j+1)} \label{gamma}%
\end{equation}
when $j\le k\le n$.
Originally, we found $\gamma_{n,k}$ using Lemma \ref{nodescoeff}, but
 an anonymous referee tells us that equation (\ref{gamma})
 is well known and can be obtained by blossoming (cf. \cite{Ra} for
an introduction to blossoms). And the Editor adds that the coefficients
 can also 
be determined via the dual functionals for B-splines. Apparently, though,
the coordinates $\gamma_{n,k}$ of $x^j$ with respect to the Bernstein basis
had not previously been  used as we do here, to define a Bernstein operator fixing $1$ and $x^{j}$. 

Next we prove convergence. 
For $l=1,\dots, j-1$, the inequalities
\[
\frac{k-j+1}{n}<\frac{k-l}{n-l}<\frac{k}{n}%
\]
can be checked by simplifying and inspection. It follows that $((k-j+1)/n)^{j}%
<t_{n,k}^{j}=\gamma_{n,k}<(k/n)^{j}$, or equivalently, that $0<k/n-t_{n,k}%
<(j-1)/n$. Thus, $B_{n,0,j}x^{m}$ converges uniformly to $x^{m}$ for
$m=0,1,2$, and by Korovkin's Theorem, $B_{n,0,j}f\rightarrow f$ uniformly for
all $f\in C[0,1].$
\end{proof}

\section{Generalized convexity}

Let $B_{n}$ denote the classical Bernstein operator defined in (\ref{defBP}%
). W.B. Temple showed in \cite{Temp54} that for a convex function $f$ the
following monotonicity property 
\begin{equation}
B_{n}f\left( x\right) \geq B_{n+1}f\left( x\right)  \label{eqSchAr}
\end{equation}
holds for all $x\in\left[ 0,1\right] .$ In \cite{Aram57} O. Aram\u{a} proved
that 
\begin{equation*}
B_{n}f\left( x\right) -B_{n+1}f\left( x\right) =\frac{x\left( 1-x\right) }{%
n\left( n+1\right) }\sum_{k=0}^{n-1}\left[ \frac{k}{n},\frac{k+1}{n+1},\frac{%
k+1}{n}\right] f \cdot\binom{n-1}{k}x^{k}\left( 1-x\right)
^{n-1-k}
\end{equation*}
where $\left[ \frac{k}{n},\frac{k+1}{n+1},\frac{k+1}{n}\right] 
f$ is the divided difference of second order, thus providing a
simple proof of Temple's result. A similar formula (see Theorem 7.5 in \cite
{Karl68}) is due to Averbach. We obtain analogous results for the
generalized Bernstein operators considered here. These generalized Bernstein
operators $B_{n}$ fix $f_{0}$ and $f_{1}$ instead of $1$ and $x$, so rather
than $(1,x)$-convexity, which is equivalent to standard convexity, the
adequate notion for our purposes is $(f_{0}, f_{1})$-convexity, to be
defined next. We shall see that for $(f_{0}, f_{1})$-convex functions $f$,
the following holds: $B_{n} f \ge B_{n+1} f \ge f$. 

Next we specialize the definition given in
\cite[p. 280]{Karl68} for an arbitrary number $n$ of functions to 
the case $n=2$. Later we will specialize it to the case $n = 1$.

\begin{definition} Let $E\subset\mathbb{R}$. A function $f: E
\rightarrow\mathbb{R}$ is called \emph{convex on }$E$\emph{\ with respect to}
a Haar pair $(f_{0},f_{1})$ if for all $x_{0},x_{1},x_{2}$ in $E$ with $%
x_{0}<x_{1}<x_{2} $, the determinant 
\begin{equation}
\operatorname{Det}_{x_{0},x_{1},x_{2}}\left( f\right) :=\det\left( 
\begin{array}{ccc}
f_{0}\left( x_{0}\right) & f_{0}\left( x_{1}\right) & f_{0}\left(
x_{2}\right) \\ 
f_{1}\left( x_{0}\right) & f_{1}\left( x_{1}\right) & f_{1}\left(
x_{2}\right) \\ 
f\left( x_{0}\right) & f\left( x_{1}\right) & f\left( x_{2}\right)
\end{array}
\right)  \label{defD}
\end{equation}
is non-negative. We shall also use the shorter expression ``$(f_{0},f_{1})$%
-convex''. Likewise, we say that $f$ is $(f_{0},f_{1})$-concave if $-f$ is $%
(f_{0},f_{1})$-convex, and $(f_{0},f_{1})$-affine if $f\in\langle
f_{0},f_{1}\rangle$.
\end{definition}

\begin{remark}
\label{nonstrict} {\rm Note that the condition $\operatorname{Det}_{x_{0},x_{1},x_{2}}%
\left( f\right) \ge0$ for all $x_{0},x_{1},x_{2}$ in $E$ with $%
x_{0}<x_{1}<x_{2}$ is equivalent to the same requirement but with $x_{0}\le
x_{1} \le x_{2}$. Of course, in the degenerate case $x_{i} = x_{i +1}$ the
determinant is zero, so it makes no difference whether or not this
possibility is included in the definition. In other words, only the ordering
of the points $x_{0},x_{1},x_{2}$ actually matters. }
\end{remark}

One of the standard definitions of convexity stipulates that the graph of $f$
must lie below the segment joining any two given points on the graph. It is
well known that an analogous characterization holds for $(f_{0},f_{1})$%
-convex functions, but with affine functions being replaced by $(f_{0},
f_{1})$-affine functions. More precisely,

\begin{proposition} Let $(f_{0},
f_{1})$ be a Haar pair.
\label{chord} Denote by $\psi_{x_{0},x_{2}}^{f}$ the unique function in $%
U_{1}:=\langle f_{0},f_{1}\rangle$ that interpolates $f$ at the points $%
x_{0}<x_{2}$, $x_{0},x_{2}\in\lbrack a,b]$, i.e., $\psi_{x_{0},x_{2}}^{f}%
\left( x_{0}\right) =f\left( x_{0}\right) $ and $\psi_{x_{0},x_{2}}^{f}%
\left( x_{2}\right) =f\left( x_{2}\right) $. Then $f$ is $(f_{0},f_{1})$%
-convex if and only if for all $x_{0},x,x_{2}$ such that $a\leq
x_{0}<x<x_{2}\leq b$, 
\begin{equation}
f\left( x\right) \leq\psi_{x_{0},x_{2}}^{f}\left( x\right) ,
\label{01convex}
\end{equation}
and in this case, for all $y\in\lbrack a,b]\setminus\lbrack x_{0},x_{2}]$, 
\begin{equation}
f\left( y\right) \geq\psi_{x_{0},x_{2}}^{f}\left( y\right) .
\label{02convex}
\end{equation}
\end{proposition}

\begin{proof}
Let $\operatorname{Det}_{x_{0},x,x_{2}}\left(  f\right)  $ be as defined in (\ref{defD}).
Observe that
\begin{align*}
\operatorname{Det}_{x_{0},x,x_{2}}\left(  f\right)   &  =\operatorname{Det}_{x_{0},x,x_{2}}\left(  f-\psi
_{x_{0},x_{2}}^{f}\right) \\
&  =-\left(  f\left(  x\right)  -\psi_{x_{0},x_{2}}^{f}\left(  x\right)
\right)  \left(  f_{1}\left(  x_{2}\right)  f_{0}\left(  x_{0}\right)
-f_{0}\left(  x_{2}\right)  f_{1}\left(  x_{0}\right)  \right)  .
\end{align*}
Since $f_{1}/f_{0}$ is strictly increasing, $f_{1}\left(  x_{2}\right)
f_{0}\left(  x_{0}\right)  -f_{0}\left(  x_{2}\right)  f_{1}\left(
x_{0}\right)  >0,$ so $\operatorname{Det}_{x_{0},x,x_{2}}\left(  f\right)  \geq0$ is equivalent
to $f\left(  x\right)  \leq\psi_{x_{0},x_{2}}^{f}\left(  x\right)  $.

Next, assume that (\ref{01convex}) holds for all $x_{0}, x, x_{2}$ such that
$a\leq x_{0}<x<x_{2}\leq b$. Suppose that for some $u\in[a,b]\setminus[x_{0},
x_{2}]$ we have $f(u) < \psi_{x_{0},x_{2}}^{f}\left(  u \right)  $. Without
loss of generality we may assume that $x_{2} < u$. We interpolate between
$x_{0}$ and $u$ to obtain a contradiction: $\psi_{x_{0},x_{2}}^{f}\left(
x_{0} \right)  = \psi_{x_{0},u}^{f}\left(  x_{0} \right)  $, while
$\psi_{x_{0},x_{2}}^{f}\left(  u \right)  > \psi_{x_{0},u}^{f}\left(  u
\right)  $. Since $\psi_{x_{0},x_{2}}^{f} - \psi_{x_{0}, u}^{f}$ has exactly
one zero, it follows that $\psi_{x_{0},x_{2}}^{f} > \psi_{x_{0}, u}^{f}$ on
$(x_{0}, b]$. In particular, $f (x_{2}) = \psi_{x_{0},x_{2}}^{f}\left(  x_{2}
\right)  > \psi_{x_{0},u}^{f}\left(  x_{2} \right)  $, which is impossible by
(\ref{01convex}) applied to $f$ and $\psi_{x_{0},u}^{f}$.
\end{proof}

Note that by inequality (\ref{01convex}), convexity is the same as $(1,x)$%
-convexity. The next result generalizes to $(f_{0},f_{1})$-convex functions,
the familiar inequality $B_{n}f\geq f$ for the classical Bernstein operator
acting on convex functions. Here it is not assumed that $B_{n}$ is defined
via an increasing sequence of nodes; it is enough to know that $%
t_{n,k}\in\lbrack a,b]$.

\begin{theorem}
\label{major} Assume that for some $n\ge1$, there is a Bernstein operator $%
B_{n}$ fixing $f_{0}$ and $f_{1}$. Then for every $(f_{0},f_{1})$-convex
function $f\in C[0,1]$ we have $B_{n}f\geq f$.
\end{theorem}

\begin{proof}
Suppose $B_{n}$ exists for some $n\geq1$, and let $\varepsilon>0$. We show
that for an arbitrary $x\in\lbrack a,b]$, $B_{n}f(x)\geq f(x)-\varepsilon$.
Assume that $x\in(a,b)$ (the cases $x=a$ and $x=b$ can be proven via obvious
changes in the notation, or just by using continuity). First, select
$\delta>0$ such that $B_{n}\delta<\varepsilon$. Next, by continuity of $f,$
choose $h>0$ so small that $[x-h,x+h]\subset\lbrack a,b]$ and $\psi
_{x-h,x+h}^{f}<f+\delta$ on $[x-h,x+h]$. Then $\psi_{x-h,x+h}^{f}<f+\delta$ on
$[a,b]$ by (\ref{02convex}), so
\[
B_{n}f(x)>B_{n}\left(  \psi_{x-h,x+h}^{f}-\delta\right)  (x)=B_{n}%
\psi_{x-h,x+h}^{f}(x)-B_{n}\delta(x)>\psi_{x-h,x+h}^{f}(x)-\varepsilon\geq
f(x)-\varepsilon,
\]
where for the last inequality we have used (\ref{01convex}).
\end{proof}

We shall use below the following characterization of $(f_{0}, f_{1})$%
-convexity, due to M. Bessen\-yei and Z. P\'ales (cf. Theorem 5, p. 388 of 
\cite{BePa}). While the result is stated there for open intervals, it also
holds for compact intervals. Note that what we call here a Haar pair (i.e., $f_0$ is strictly positive and 
$f_{1}/f_{0}$ strictly increasing) is called in \cite{BePa} a positive regular pair.

\begin{theorem}
\label{bepa} Let $(f_{0},
f_{1})$ be a Haar pair and let $I:= (f_{1}/f_{0})([a,b])$. Then $f\in C[a,b]$ is $(f_{0},
f_{1})$-convex if and only if $(f/f_{0}) \circ(f_{1}/f_{0})^{-1} \in C\left(
I\right) $ is convex in the standard sense.
\end{theorem}

\begin{example}
{\rm Consider the Bernstein operator $B_{n,0,j}$ from Proposition \ref{0j},
defined on $C[0,1]$ and fixing 1 and $x^{j}$. It is easy to see from Theorem 
\ref{bepa} that for $s\in(0,j)$, the function $x^{s}$ is $(1,x^{j})$%
-concave, while if $s\in(j,\infty)$, $x^{s}$ is $(1,x^{j})$-convex.
Therefore, by Theorem \ref{major}, for all $x\in[0,1]$ we have $B_{n,0,j}
x^{s} \le x^{s}$ if $s\in(0,j)$ and $B_{n,0,j} x^{s} \ge x^{s}$ when $%
s\in(j,\infty)$.}
\end{example}

Our next objective is to obtain an analog of Aram\u{a}'s result (presented
at the beginning of this section) for generalized Bernstein operators 
$B_{n}$. Here the interlacing property of nodes is used in an essential way.

\begin{proposition}
\label{Prop0} Let $s_{k}, s_{k+1}, s_{k+2}\in[a,b]$ be such that $%
s_{k}<s_{k+1}<s_{k+2}$, and assume that $G_{k}:C\left[ c,d\right] \rightarrow%
\mathbb{R}$ is a functional of the form 
\begin{equation*}
G_{k}\left( f\right) =a_{k}f\left( s_{k}\right) +b_{k}f\left( s_{k+1}\right)
+c_{k}f\left( s_{k+2}\right) ,
\end{equation*}
satisfying $G_{k}\left( f_{0}\right) =G_{k}\left( f_{1}\right) =0$. Then $%
b_{k}\geq0$ if an only if $G_{k}\left( f\right) \leq0$ for all $(f_{0}$, $%
f_{1})$-convex functions $f\in C\left[ t_{k},t_{k+2}\right] $.
\end{proposition}

\begin{proof}
Let $f$ be $(f_{0}$, $f_{1})$-convex, and let $\psi_{s_{k},s_{k+2}}^{f}$ be
the function in $\left\langle f_{0},f_{1}\right\rangle $ that interpolates $f$
at the points $s_{k}$ and $s_{k+2}.$ By (\ref{01convex}),
\[
G_{k}\left(  f\right)  =G_{k}\left(  f-\psi_{s_{k},s_{k+2}}^{f}\right)
=b_{k}\left(  f -\psi_{s_{k},s_{k+2}}^{f}\right)  \left(  s_{k+1}\right)
\leq0
\]
if and only if $b_{k} \ge0$.
\end{proof}

\begin{theorem}
\label{decreasing} Under the same hypotheses and with the same notation as
in Theorem \ref{ThmExt}, let the linear functionals $G_{k},k=1,...,n-1$, be
defined by 
\begin{equation*}
G_{k}\left( f\right) =f\left( t_{n,k}\right) \alpha_{n,k}\frac
{p_{n,k}^{\left( k\right) }\left( a\right) }{p_{n+1,k}^{\left( k\right)
}\left( a\right) }-f\left( {t}_{n+1,k}\right) {\alpha}_{n+1,k}+f\left( t_{n,
k-1}\right) \alpha_{n,k-1}\frac{p_{n,k-1}^{\left( n+1-k\right) }\left(
b\right) }{p_{n+1,k}^{\left( n+1-k\right) }\left( b\right) }.
\end{equation*}
Then 
\begin{equation*}
B_{n}f-B_{n+1}f=\sum_{k=1}^{n}G_{k}\left( f\right) \cdot p_{n+1,k}.
\end{equation*}
In particular, if $f$ is $\left( f_{0},f_{1}\right) $-convex then $%
B_{n}f-B_{n+1}f\geq0$.
\end{theorem}

\begin{proof}
Recall that
\[
B_{n}f=\sum_{k=0}^{n}f\left(  t_{n,k}\right)  \alpha_{n,k}p_{n,k}\text{ and
}B_{n+1}f=\sum_{k=0}^{n+1}f\left(  {t}_{n+1,k}\right)  {\alpha}_{n+1,k}%
p_{n+1,k}%
\]
where $t_{n,0}={t}_{n+1,0}=a$, $t_{n,n}={t}_{n+1,n+1}=b$, and $t_{n,k-1}%
<{t}_{n+1,k}<t_{n,k}$ for $k=1,...,n$. Using Proposition \ref{Prop1} we
obtain
\begin{align*}
B_{n}f-B_{n+1}f  &  =\sum_{k=0}^{n}f\left(  t_{n,k}\right)  \alpha_{n,k}%
\frac{p_{n,k}^{\left(  k\right)  }\left(  a\right)  }{p_{n+1,k}^{\left(
k\right)  }\left(  a\right)  }p_{n+1,k}\\
&  +\sum_{k=0}^{n}f\left(  t_{n,k}\right)  \alpha_{n,k}\frac{p_{n,k}^{\left(
n-k\right)  }\left(  b\right)  }{p_{n+1,k+1}^{\left(  n-k\right)  }\left(
b\right)  }p_{n+1,k+1}-\sum_{k=0}^{n+1}f\left(  {t}_{n+1,k}\right)  {\alpha
}_{n+1,k}p_{n+1,k}.
\end{align*}
It follows from (\ref{coeff0n}) that the first summands (corresponding to
$k=0$) of the first and the last sum are the same, so they cancel out.
Likewise, the $n$-th summand of the second sum and the $\left(  n+1\right)
$-st summand of the last sum cancel out. Thus
\[
B_{n}f-B_{n+1}f=
\]%
\[
\sum_{k=1}^{n}p_{n+1,k}\left[  f\left(  t_{n,k}\right)  \alpha_{n,k}%
\frac{p_{n,k}^{\left(  k\right)  }\left(  a\right)  }{p_{n+1,k}^{\left(
k\right)  }\left(  a\right)  }-f\left(  {t}_{n+1,k}\right)  {\alpha}%
_{n+1,k}+f\left(  t_{n,k-1}\right)  \alpha_{n,k-1}\frac{p_{n,k-1}^{\left(
n+1-k\right)  }\left(  b\right)  }{p_{n+1,k}^{\left(  n+1-k\right)  }\left(
b\right)  }\right].
\]
Finally, let $f$ be $\left(  f_{0},f_{1}\right)  $-convex. Taking
$s_{k}=t_{n,k-1},s_{k+1}={t}_{n+1,k}$, and $s_{k+2}=t_{n,k}$ in Proposition
\ref{Prop0}, we get $B_{n}f-B_{n+1}f\geq0$.
\end{proof}

A very natural question, not touched upon here, is under which conditions a
sequence of Bernstein operators for $(f_{0},f_{1})$ converges to the
identity. It follows from Theorems \ref{major} and \ref{decreasing} that if $%
f$ is $(f_{0},f_{1})$-convex, then the sequence $\{B_{n}f\}_{n=1}^{\infty}$
monotonically converges to some function $g\geq f$ (assuming that a sequence
of functions $f_{0},f_{1},f_{2},...$ are given such that $\left\langle
f_{0},...,f_{n}\right\rangle $ is an extended Chebyshev space of dimension $%
n+1$ for each $n\in\mathbb{N}).$ But we have not determined which conditions
will ensure that $g=f$. In this regard, we expect the strict interlacing
property of nodes to be useful, since it entails, in a qualitative sense,
that the sampling of functions is not ``too biased''.

\section{Total positivity and generalized convexity}

Let $B_{n}:C\left[ a,b\right] \rightarrow U_{n}$ be a Bernstein operator for
the pair $(f_{0},f_{1}).$ In Section 3 we proved that $B_{n}f\geq f$ for all 
$\left( f_{0},f_{1}\right) $-convex functions $f\in C\left[ a,b\right] .$
This did not require an increasing sequence of nodes; it was enough to know
that $t_{n,k}\in\lbrack a,b].$

In this section we show that $B_{n}f$ is $\left( f_{0},f_{1}\right) $-convex
for every $\left( f_{0},f_{1}\right) $-convex function $f\in C\left[ a,b%
\right] $, provided that the nodes $t_{n,0},...,t_{n,n}$ are \emph{increasing%
} and $U_{n}$ is an extended Chebyshev space over $\left[ a,b\right] ,$ and
a similar result holds for the so-called $g$-monotone functions. These
statements  follow directly from the more general results presented in \cite
{Karl68} concerning shape preserving properties of linear transformations
with totally positive kernels. The connection between total positivity and
shape preserving properties of bases is a classical subject that has been
widely described, see e.g. \cite{Karl68} or the more recent survey \cite
{CaPe96}.

The following definitions come from \cite{Karl68}. Let $X$ and $Y$ be
subsets of $\mathbb{R}.$ A function $K:X\times Y\rightarrow\mathbb{R}$ is
called \emph{sign-consistent of order }$m$ if there exists an $%
\varepsilon_{m}\in\left\{ -1,1\right\} $ such that 
\begin{equation}
\varepsilon_{m}\det\left( 
\begin{array}{cccc}
K\left( x_{1},y_{1}\right) & K\left( x_{1},y_{2}\right) & ... & K\left(
x_{1},y_{m}\right) \\ 
K\left( x_{2},y_{1}\right) & K\left( x_{2},y_{2}\right) & ... & K\left(
x_{2},y_{m}\right) \\ 
&  &  &  \\ 
K\left( x_{m},y_{1}\right) & ... & .... & K\left( x_{m},y_{m}\right)
\end{array}
\right) \geq0  \label{dettotal}
\end{equation}
for all $x_{1}<x_{2}<...<x_{m}$ in $X$ and $y_{1}<y_{2}<...<y_{m}$ in $Y.$
If $\varepsilon_{m}=1$ we shall call $K$ \emph{positive of order} $m.$ A
function $K$ is \emph{totally positive} if it is positive of all orders $m$
with $m\in\mathbb{N},m\geq1.$ Similarly, if one has strict positivity in (%
\ref{dettotal}) then $K$ is called \emph{strictly sign-consistent of order} $%
m,$ and if in addition $\varepsilon_{m}=1$, then $K$ is \emph{strictly
positive of order} $m.$ Strict total positivity means that $K$ is strictly
positive of all orders $m\in\mathbb{N},m\geq1.$

The following result is well-known, see e.g. \cite[p. 358]{Mazu05}, or the
proof presented in \cite[pp. 342--344]{CLM}:

\begin{theorem}
Let $U_{n}\subset C^{n}\left[ a,b\right] $ be an extended Chebyshev space
over $\left[ a,b\right] $ and let $p_{n,k},k=0,...,n,$ be a non-negative
Bernstein basis for $\left[ a,b\right] .$ Then $K:[a,b]\times\left\{
0,....,n\right\} \rightarrow\mathbb{R}$, defined by 
\begin{equation}
K\left( x,k\right) :=p_{n,k}\left( x\right)  \label{defK},
\end{equation}
is totally positive, and $K$ is strictly totally positive on $\left(
a,b\right) \times\left\{ 0,....,n\right\} .$
\end{theorem}

Following the notation of \cite{CMP04}, \cite{CMP07}, \cite{MPS01}, we can
deduce from the previous result that a non-negative Bernstein basis of an
extended Chebyshev system over $\left[ a,b\right] $ is totally positive on $%
\left[ a,b\right] $, i.e.,  
 a B-basis.

We cite from \cite[p. 284]{Karl68} the following result (specialized to the
case of two functions $F_{0},F_{1}$ instead of a family $F_{1},...,F_{m}$).

\begin{theorem}
\label{ThmC}Let $X$ and $Y$ be subsets of $\mathbb{R,}$ let $F_{0},F_{1}$ be
functions on $Y$ and let $K:X\times Y\rightarrow\mathbb{R}$ be continuous,
and positive of order $3$. Let $\mu$ be a non-negative sigma-finite measure
and $B_{K}:C\left( Y\right) \rightarrow C\left( X\right) $ be defined by
\begin{equation*}
B_{K}\left( F\right) \left( x\right) :=\int_{Y}K\left( x,y\right) F\left(
y\right) d\mu\left( y\right) .
\end{equation*}
If $F$ is $(F_{0},F_{1})$-convex, then $B_{K}\left( F\right) $ is $%
(B_{K}F_{0},B_{K}F_{1})$-convex.
\end{theorem}

From this we conclude:

\begin{theorem}
\label{Thmcconvex} Let $U_{n}$ be an extended Chebyshev space over $\left[
a,b\right] .$ Assume there exists a Bernstein operator $B_{n}:C\left[ a,b%
\right] \rightarrow U_{n}$ fixing $f_{0}$ and $f_{1}$, with increasing nodes 
$t_{n,0}\leq....\leq t_{n,n}.$ If $f\in C\left[ a,b\right] $ is $\left(
f_{0},f_{1}\right) $-convex, then $B_{n}\left( f\right) $ is $\left(
f_{0},f_{1}\right) $-convex.
\end{theorem}

\begin{proof}
Put $X=\left[  a,b\right]  $ and $Y:=\left\{  0,....,n\right\}  $. Define the
function $\varphi:Y\rightarrow X$ by $\varphi\left(  k\right)  :=t_{n,k}$, for
$k=0,...,n$. Observe that $\varphi$ is monotone increasing (though perhaps
not strictly), so it is order preserving. Next, set $\mu:=\sum_{k=0}^{n}%
\alpha_{n,k}\delta_{k}$, where the $\alpha_{n,k}$ are the positive
coefficients defining $B_{n}$ and $\delta_{k}$ is the Dirac measure at the
point $k\in\left\{  0,...,n\right\}  .$ With $K(x,k):=p_{n,k}(x)$, we obtain,
for every $F\in C\left(  Y\right)  $,
\[
B_{K}\left(  F\right)  \left(  x\right)  :=\int_{Y}K\left(  x,y\right)
F\left(  y\right)  d\mu\left(  y\right)  =\sum_{k=0}^{n}F\left(  k\right)
\alpha_{n,k}p_{n,k}\left(  x\right)  .
\]
Now let $f\in C\left(  X\right)  ,$ and define $F:=f\circ\varphi\in C\left(
Y\right)  .$ Then
\begin{equation}
B_{K}\left(  f\circ\varphi\right)  \left(  x\right)  =\sum_{k=0}^{n}f\left(
t_{n,k}\right)  \alpha_{n,k}p_{n,k}\left(  x\right)  =B_{n}\left(  f\right)
\left(  x\right)  . \label{eqnn}%
\end{equation}
If $f\in C\left(  X\right)  $ is $\left(  f_{0},f_{1}\right)  $-convex, then
$F=f\circ\varphi$ is $\left(  f_{0}\circ\varphi,f_{1}\circ\varphi\right)
$-convex, since $\varphi$ preserves order (cf. Remark \ref{nonstrict}).
Putting $F_{j}=f_{j}\circ\varphi$ for $j=0,1,$ an application of Theorem
\ref{ThmC} shows that $B_{K}\left(  F\right)  $ is $(B_{K}F_{0},B_{K}F_{1}%
)$-convex. By formula (\ref{eqnn}) and the property that $B_{n}$ fixes $f_{0}$
and $f_{1}$ one obtains
\[
B_{K}F_{j}=B_{n}\left(  f_{j}\right)  =f_{j}%
\]
for $j=0,1.$ Thus $B_{K}\left(  F\right)  =B_{n}\left(  f\right)  $ is
$\left(  f_{0},f_{1}\right)  $-convex.
\end{proof}

In a similar way it is possible to obtain generalized monotonicity properties of the
Bernstein operator.

\begin{definition}
\textrm{Let $g>0$. We say that $f$ is }$g$\textrm{-increasing on }$\left[ a,b%
\right] $\textrm{\ if $f/g$ is increasing on }$\left[ a,b\right] ,$ i.e.
if 
\begin{equation*}
\frac{f\left( x_{0}\right) }{g\left( x_{0}\right) }\leq\frac{f\left(
x_{1}\right) }{g\left( x_{1}\right) }
\end{equation*}
for all $x_{0}<x_{1}$ in $\left[ a,b\right] .$
\end{definition}

The notions of 
$g$-decreasing and $g$-monotone are the obvious ones.
Now $f:\left[ a,b\right] \rightarrow \mathbb{R}
$ is $g$-increasing on\textrm{\ }$\left[ a,b\right] $ if and only if for all $%
x_{0},x_{1} $ in $\left[ a,b\right] $ with $x_{0}<x_{1}$ the determinant 
\begin{equation*}
\det\left( 
\begin{array}{cc}
g\left( x_{0}\right) & g\left( x_{1}\right) \\
f\left( x_{0}\right) & f\left( x_{1}\right) 
\end{array}
\right)
\end{equation*}
is non-negative. 
But this condition is just the definition of convexity with respect to
one function (the function $g$), obtained by specializing to $n=1$ the definition given in
\cite[p. 280]{Karl68}. Specializing also Theorem 3.3 in \cite[p. 284]{Karl68} to one function (cf. Theorem \ref{ThmC} above
in the case of two functions) and arguing as in  the
proof of Theorem \ref{Thmcconvex}, we obtain:

\begin{theorem}
\label{Thmmon}Let $U_{n}$ be an extended Chebyshev space over $\left[ a,b%
\right] .$ Assume that there exists a Bernstein operator $B_{n}:C\left[ a,b\right]
\rightarrow U_{n}$ with
increasing nodes $t_{n,0}\leq ....\leq t_{n,n}$, fixing
 $f_{0}$ and $f_{1}$.  If $g\in \left\langle
f_{0},f_{1}\right\rangle $ is strictly positive and $f\in C\left[ a,b\right] $ is $%
g$-monotone, then $B_{n}\left( f\right) $ is $g$-monotone.
\end{theorem}

One of the referees points out that  Theorems \ref
{Thmmon} and \ref{Thmcconvex} can be proven without referring to Theorem 3.3 in \cite[p.
284]{Karl68},  using instead elementary shape preserving properties of totally
positive bases, as described in the surveys \cite{CaPe96} and \cite{Good96}.
Moreover, Theorems \ref{Thmmon} and \ref{Thmcconvex} in the special case $%
f_{0}=1$ and $f_{1}\left( x\right) =x$ follow immediately from   \cite[Corollaries
3.7 and 3.8, p. 162]{Good96}.

\section{Normalized Bernstein bases and existence of Bernstein operators}

Let $b_{n,k},k=0,...,n$, be a basis of a given subspace $U_{n}\subset C\left[
a,b\right] $ of dimension $n+1.$ The basis $b_{n,k},k=0,...,n$ is \emph{%
totally positive} if the kernel $K\left( x,k\right) :=b_{n,k}\left( x\right) 
$ is totally positive (in particular, the functions $b_{n,k}$ are
non-negative). Suppose now that $U_{n}$ contains the constant function $1.$
Then the basis $b_{n,k},k=0,...,n$ is called \emph{normalized} if 
\begin{equation*}
1=\sum_{k=0}^{n}b_{n,k}\left( x\right)
\end{equation*}
for all $x\in \left[ a,b\right] .$ Normalized totally positive bases are
important in geometric design due to their good shape preserving properties.
J.-M. Carnicer and J.-M. Pe\~{n}a have shown that  normalized, totally
positive Bernstein bases are optimal from the shape preservation viewpoint, see \cite{CaPe93}, \cite{CaPe94}, \cite
{Good96}. Moreover it was proven in \cite{CMP04}, and independently in \cite
{Mazu05}, that a subspace $U_{n}$ of $C^{n}\left[ a,b\right] $ of dimension $%
n+1$ containing the constant functions, possesses a normalized, totally
positive Bernstein basis, provided that {\em both} $U_{n}$ and the space of all
derivatives $U_{n}^{\prime }:=\left\{ f^{\prime }:f\in U_{n}\right\} $ are
extended Chebyshev spaces over $\left[ a,b\right] .$

From this point of view it is natural to conjecture that to define
a well-behaved Bernstein operator (with increasing nodes)
fixing $f_0$ and $f_1$, it is enough to assume  that $U_{n}$ possesses a normalized totally positive Bernstein
basis and that $\left\langle f_{0}\right\rangle ,\left\langle
f_{0},f_{1}\right\rangle $ are extended Chebyshev systems. However, we shall
show by a counterexample that this is not true. We refer to \cite{AKR08} for
a more detailed discussion under which conditions there might exist a
Bernstein operator fixing a pair $f_{0},f_{1}\in U_{n}.$

Consider the linear space $\mathcal{U}_3:= \langle 1,x,\cos x,\sin x
\rangle$ 
over the interval $\left[ 0,b\right]$. This space has been previously
studied  
by several authors, see the references in \cite{CMP04} or \cite{MPS01}. 
It is well known
that $\mathcal{U}_3$ and  $\mathcal{U}_3^{\prime }$ (the space of all derivatives 
of functions in $\mathcal{U}_3$) are extended
Chebyshev spaces over $\left[ 0,b\right] $ for every $b\in \left( 0,2\pi
\right) .$ Thus $\mathcal{U}_3$ possesses a normalized totally positive Bernstein
basis for every $b\in \left( 0,2\pi \right) .$ By \cite{CMP04} this entails
that the critical length of $\mathcal{U}_3$ for design purposes is $2\pi $. However,
we show in Theorem \ref{ThmExamp5} that for $b$ sufficiently close to $2\pi $
(say, $b\geq 4.5$) there is no Bernstein operator from $C[0,b]$ to $\mathcal{U}_3$ fixing $1$ and $x$. 

The obstruction for employing Corollary \ref{chain} is due to the fact that
neither $\langle 1,x,\cos x\rangle $ nor $\langle 1,x,\sin x\rangle $ are
extended Chebyshev spaces over $[0,b]$ for all $b<2\pi $ (for instance, $%
\sin x-x$ has a zero of order 3 at 0) so the chain of nested spaces cannot
be continued beyond $U_{1}=\left\langle 1,x\right\rangle $. By Corollary \ref
{space}, it is nevertheless possible to construct a Bernstein operator
fixing {\em some} Haar pair of functions $g_{0},g_{1}\in \mathcal{U}_3^{\prime }=\left\langle
1,\cos x,\sin x\right\rangle $. Hence, by Theorem \ref{ThmExt} there is a corresponding
Bernstein operator from $C[0,b]$ to $\mathcal{U}_3$, fixing $g_{0}$ and $g_{1}$, with strictly
interlacing nodes.

\begin{theorem}
\label{ThmExamp5}
Given $b\in (0, 2\pi)$, let $\rho _{0}$ be the first positive root of $b\mapsto \sin b-b \cos b,$
($\rho _{0}\approx 4.4934$). Let $\mathcal{U}_3=\left\langle 1,x,\cos x,\sin x\right\rangle $, 
 $f_{0}=1$, and  $f_{1}\left( x\right) =x.$ Then for every $b\in \left(
0,\rho _{0}\right] $ there exists a Bernstein operator 
$B_3 : C[0,b] \to \mathcal{U}_3$ fixing 
 $1$ and $x$. The nodes $t_{0}\left(
b\right) ,t_{1}\left( b\right) ,t_{2}\left( b\right)$, and  
$t_{3}\left( b\right) $
satisfy the following inequalities: 
\begin{eqnarray*}
0 &=&t_{0}\left( b\right) <t_{1}\left( b\right) <t_{2}\left( b\right)
<t_{3}\left( b\right) =b\text{ for }b\in \left( 0,\pi \right)  \\
0 &=&t_{0}\left( b\right) <t_{1}\left( b\right) = \pi /2 =t_{2}\left( b\right)
<t_{3}\left( b\right) =b\text{ for }b=\pi  \\
0 &=&t_{0}\left( b\right) <t_{2}\left( b\right) <t_{1}\left( b\right)
<t_{3}\left( b\right) =b\text{ for }b\in \left( \pi ,\rho _{0}\right) .
\end{eqnarray*}
 The operator $B_3$ preserves
convex functions whenever $b\in \left( 0,\pi \right] ,$ but not for $b\in
\left( \pi ,\rho _{0}\right].$ Finally, if $b\in \left( \rho _{0},2\pi \right) $, then there does not exist a Bernstein
operator fixing 
 $1$ and $x$.
\end{theorem}

To keep computations simple we shall present  first two general
propositions and a definition: We say that a subspace $%
U_{n}\subset C\left[ a,b\right] $ is \emph{symmetric} if $f\in U_{n}$
implies that the function $F$ defined by $F\left( x\right) :=f\left(
a+b-x\right) $ is in $U_{n}.$

Assume that $U_{n}\subset C^{n}\left[ a,b\right] $ is a symmetric, extended
Chebyshev space over $\left[ a,b\right] .$ Let $p_{n,k},k=0,....,n$ be a
non-negative Bernstein basis of $U_{n}$, and let $\beta_{0},...,\beta_{n}$
and $\gamma_{0},...,\gamma_{n}$ be constants such that $1=\sum_{k=0}^{n}%
\beta _{k}p_{n,k}\left( x\right) $ and $x=\sum_{k=0}^{n}\gamma_{k}p_{n,k}%
\left( x\right) $ for all $x\in\left[ a,b\right] .$ 

\begin{proposition}
\label{Proptt}Suppose that $U_{n}\subset C^{n}\left[ a,b\right] $ is a
symmetric, extended Chebyshev space over $\left[ a,b\right] $ containing the
constant function $f_{0}=1$ and the identity function $f_{1}\left( x\right)
=x.$ If there exists a Bernstein operator $B_{n}$ fixing $f_{0}$ and $f_{1}$%
, then the following equalities hold for the coefficients $\beta_{k}$ and
the nodes $t_{n,k}$, whenever $k=0,...,n$: 
\begin{equation}
\beta_{k} = \beta_{n-k}\mbox{ \ \ \ and \ \ \ }t_{n,k}+t_{n,n-k}=a+b.
\label{eqsym}
\end{equation}
\end{proposition}

\begin{proof}
Let  $\tilde U_n:= \{f(a + b -x) : f\in U_n\}$. Given $g\in \tilde U_n$, define $\tilde B_n g(x) := B_n g (a + b - x)$. Since by symmetry
$\tilde U_n = U_n$, and both operators $\tilde B_n$ and $B_n$ 
fix the affine functions on $[a,b]$,  by uniqueness of the Bernstein
operator we have  
$\tilde B_n  =B_n$. Thus, $B_n f(x) = B_n f (a + b - x)$ for every
$f\in U_n$. Likewise, let $p_{n,k},k=0,....,n$ be a non-negative Bernstein basis of $U_{n}$. By
suitably rescaling we may impose the following normalization: 
\begin{equation}
p_{n,k}\left(  \frac{a+b}{2}\right)  =1 \label{eqpk1}%
\end{equation}
for all $k=0,....,n$. Counting zeros at the endpoints, and using uniqueness of the Bernstein basis (up to a normalizing constant), we conclude that
\begin{equation}
p_{n,k}\left(  a+b-x\right)  =p_{n,n-k}\left(  x\right)  . \label{eqsymmet}
\end{equation}
And now  (\ref{eqsym}) follows from 
 (\ref{eqsymmet}) and the fact that $B_n f(x) = B_n f (a + b - x)$ for every
 $f\in U_n$. 
\end{proof}

\begin{proposition}
\label{PropTr} Let $p_{n,k},k=0,...,n,$ be a Bernstein basis, and for $f\in
U_{n}$ let $\beta_{0},...,\beta_{n}$ be the coefficients in the expression 
\begin{equation}
f=\sum_{k=0}^{n}\beta_{k}p_{n,k}.  \label{eqbneu}
\end{equation}
Then $p_{n,n}\left( b\right) \beta_{n}=f\left( b\right) $ and 
\begin{equation}
f^{\prime}\left( b\right) =\beta_{n-1}p_{n,n-1}^{\prime}\left( b\right)
+\beta_{n}p_{n,n}^{\prime}\left( b\right).  \label{eqfstrich}
\end{equation}
\end{proposition}

\begin{proof}
For the first statement insert $x=b$ in (\ref{eqbneu}), for the second take
the derivative of $f$ in (\ref{eqbneu}) and then insert $x=b.$
\end{proof}

Now we turn to the proof of Theorem \ref{ThmExamp5}:

\begin{proof}
We define a Bernstein basis for $\mathcal{U}_3$ over  $\left[ 0,b\right] $ by 
setting
\begin{eqnarray*}
p_{3,3}\left( x\right) &=&x-\sin x, \\
p_{3,2}\left( x\right) &=&\left( b-\sin b\right) \left( 1-\cos x\right)
-\left( 1-\cos b\right) \left( x-\sin x\right) , \\
p_{3,1}\left( x\right) &=&p_{3,2}\left( b-x\right) \text{ and }p_{3,0}\left(
x\right) =p_{3,3}\left( b-x\right) .
\end{eqnarray*}

Let us check that this is indeed a Bernstein basis.
We claim that the four given functions are positive 
on $(0, \pi)$ and have the
required number of zeros at the endpoints. This is clear
for $p_{3,3}$, and hence for $p_{3,0}$. Regarding $p_{3,2}$, 
note that it has
a zero of order 2 at 0 and another zero at $b$. Since  $\left\langle 1,x,\cos x,\sin x\right\rangle$ is an extended Chebyshev space whenever
$b\in (0, 2\pi)$, 
 $p_{3,2}$ has at most 3 zeros in $[0,b]$. So it
has no zeros in $(0,b)$  and the zero at $b$ is of order
exactly  one. Since $p^{\prime\prime}_{3,2} (0) > 0$, $p_{3,2} \ge 0$ on $[0, b]$. Thus, $p_{3,1} \ge 0$ on $[0, b]$.

Recall that Bernstein bases are unique up to multiplicative constants.
So to prove that a Bernstein operator does not exist, it is sufficient
to consider the preceding basis. On the other hand, to prove that a Bernstein
operator does exist, we need to exhibit nodes $t_k$ in $[0,b]$ and positive
coefficients $\alpha_k$, for $k = 0,1,2,3$.
Now let $1=\sum_{k=0}^{3}\beta _{k}p_{3,k}$ and $x=\sum_{k=0}^{3}\gamma
_{k}p_{3,k}.$ Lemma \ref{nodescoeff} tells us what the nodes and coefficients
must be if $B_3$ exists. 
By Proposition \ref{PropTr} we have 
\begin{equation*}
\beta _{3}=\frac{1}{b-\sin b}\text{ and }\gamma _{3}=\frac{b}{b-\sin b},
\end{equation*}
so $t_3(b) := \gamma _{3}/\beta _{3}=b$. 
It follows from (\ref{eqfstrich}) that 
\begin{equation}\label{beta2}
\beta _{2}p_{3,2}^{\prime }\left( b\right) =-\beta _{3}p_{3,3}^{\prime
}\left( b\right) \text{ and }\gamma _{2}p_{3,2}^{\prime }\left( b\right)
=1-\gamma _{3}p_{3,3}^{\prime }\left( b\right) .
\end{equation}
Thus   
\begin{equation*}
t_{2}\left( b\right) 
:= 
\frac{\gamma _{2}}{\beta _{2}}
=
\frac{-\gamma
_{3}p_{3,3}^{\prime }\left( b\right) }{-\beta _{3}p_{3,3}^{\prime }\left(
b\right) } + \frac{1}{-\beta _{3}p_{3,3}^{\prime }\left(
b\right) }
=
b-\frac{b-\sin b}{1-\cos b}.
\end{equation*}
Now (\ref{eqsym}) implies that 
\begin{equation*}
t_{1}\left( b\right) =\frac{b-\sin b}{1-\cos b}.
\end{equation*}
We see that $t_{2}\left( b\right) - t_{1}\left( b\right) > 0$ if and only if 
$
g(b):= 2\sin b -b \cos b - b >0.
$
Elementary calculus shows that $g > 0$ on $(0, \pi)$, $g(\pi) = 0$, and
$g < 0$ (at least) on $(\pi, 3 \pi/2)$. If $b=\pi $, then  $t_{1}\left( \pi\right) =t_{2}\left( \pi\right) =\pi /2.$
Furthermore, $t_{2}\left( b\right) <0$ whenever $\sin b-b\cos b<0,$ so by
Lemma \ref{nodescoeff}, for $b\in \left(
\rho _{0},2\pi \right) $ there does not exists a Bernstein operator.
To see that such operator exists when $b\in (0,\rho_0)$, note that
since $f_0 = 1$, by (\ref{coeff}) we have $\alpha_k = \beta_k$, so it is
enough to show that  $\beta_k > 0$ for $k = 0,1,2,3$. Since  $\beta_0 = \beta_3$ and  $\beta_1 = \beta_2$ by  Proposition \ref{Proptt}, and $\beta_3 > 0$, it suffices to prove that
$\beta_2 > 0$. Now
from equation (\ref{beta2}) we get 
\begin{equation*}
\beta_{2} = 
- \left(\frac{1 - \cos b}{b - \sin b}\right)
\frac{1}{ b \sin b - 2 + 2 \cos b },
\end{equation*}
so $\beta_2 > 0$ if and only if $ b \sin b - 2 + 2 \cos b < 0$. Elementary
calculus shows that this is the case for every $b\in (0, 2 \pi)$,
and in particular, for every $b\in (0, \rho_0]$.

Regarding the convexity assertions, if $b\in \left( 0,\pi \right]$ then the Bernstein operator $B_{3}$
preserves convexity by Theorem \ref{Thmcconvex}. Next, fix $b\in \left( \pi
,\rho _{0}\right]$,   write $t_{1}=t_{1}\left( b\right) $, 
$t_{2}=t_{2}\left( b\right)$, and consider the convex function $f\left(
x\right) =\left( x-t_{1}\right) \left( x-t_{2}\right)$. Since 
$t_{1}+t_{2}=b$, we have 
$
f\left( 0\right) =f\left( b\right) =t_{1}\left( b-t_{1}\right).
$
By Proposition \ref{Proptt}, $\beta _{0}=\beta _{3},$ so 
\begin{equation*}
B_{3}f\left( x\right) =f\left( 0\right) \beta _{0}p_{3,0}\left( x\right)
+f\left( b\right) \beta _{3}p_{3,3}\left( x\right) =\beta _{0}f\left(
0\right) \left( p_{3,0}\left( x\right) +p_{3,3}\left( x\right) \right) .
\end{equation*}
Using $\beta _{0}f\left( 0\right) > 0$  we see that $B_{3}f$ is
convex if and only if $F:=p_{3,0}+p_{3,3}$ is convex. A direct computation
shows that $F\left( x\right) =b-\sin \left( b-x\right) -\sin x,$ so 
\begin{equation*}
F^{\prime \prime }\left( x\right) =\sin \left( b-x\right) +\sin x.
\end{equation*}
Thus $F^{\prime \prime }\left( 0\right) =\sin b<0$, since  $b\in
\left( \pi ,2\pi \right) .$ By continuity, $F^{\prime \prime }\left(
x\right) <0$ for all $x$ in a small neighborhood of $0,$ so $F$ is not
convex. 
\end{proof}

 We see that  non-increasing nodes can lead to pathological behavior on the part of the generalized Bernstein operator defined by them. Thus,
either additional conditions are imposed to avoid this situation (in this
case, a smaller value of $b$) or else the requirement
that the operator fix $1$ and $x$ must be given up. 
 On the other hand, there are natural examples where  nodes are 
increasing, but not strictly. In addition to Proposition  \ref{0j}, where the
node at zero is repeated $j$ times, we mention the case $b=\pi$ in Theorem \ref{ThmExamp5}.  There, the second and third nodes are equal 
(to $\pi/2$). It is shown in \cite{AR08} that the generalized Bernstein
operator on $\mathcal{U}_3$ approximates some convex functions (such as
$| x - \pi/2|$) on $[0,\pi]$ better than the
standard polynomial Bernstein operator on $\langle 1, x, x^2, x^3\rangle$ (while the latter operator approximates some other functions better).

\end{document}